\documentclass[12pt]{amsart}
\usepackage{a4wide}
\usepackage[utf8]{inputenc}
\usepackage{amsmath}
\usepackage{amsfonts}
\usepackage{amssymb}
\usepackage{amsthm}
\usepackage{subfigure}
\usepackage{color}

\newtheorem{theo}{Theorem}[section]
\newtheorem{prop}[theo]{Proposition}
\newtheorem{coro}[theo]{Corollary}
\newtheorem{lemm}[theo]{Lemma}

\theoremstyle{definition}

\theoremstyle{remark}
\newtheorem{rema}[theo]{Remark}
\newcommand{\Op}{\operatorname{Op}}

\newcommand{\nwc}{\newcommand}
\nwc{\eps}{\epsilon}
\nwc{\ep}{\epsilon}
\nwc{\vareps}{\varepsilon}
\nwc{\Oph}{\operatorname{Op}_\hbar}
\nwc{\la}{\langle}
\nwc{\ra}{\rangle}

\nwc{\mf}{\mathbf} 
\nwc{\blds}{\boldsymbol} 
\nwc{\ml}{\mathcal} 

\nwc{\defeq}{\stackrel{\rm{def}}{=}}

\nwc{\cE}{\ml{E}}
\nwc{\cN}{\ml{N}}
\nwc{\cO}{\ml{O}}
\nwc{\cP}{\ml{P}}
\nwc{\cU}{\ml{U}}
\nwc{\cV}{\ml{V}}
\nwc{\cW}{\ml{W}}
\nwc{\tU}{\widetilde{U}}
\nwc{\IN}{\mathbb{N}}
\nwc{\IQ}{\mathbb{Q}}
\nwc{\IR}{\mathbb{R}}
\nwc{\IZ}{\mathbb{Z}}
\nwc{\IC}{\mathbb{C}}
\nwc{\IT}{\mathbb{T}}
\nwc{\tP}{\widetilde{P}}
\nwc{\tPi}{\widetilde{\Pi}}
\nwc{\tV}{\widetilde{V}}
\nwc{\supp}{\operatorname{supp}}
\nwc{\rest}{\restriction}

\begin{document}

\title[Two-microlocal regularity of quasimodes on $\IT^2$]{Two-microlocal regularity of quasimodes on the torus}

\author[Fabricio Maci\`a]{Fabricio Maci\`a}
\author[Gabriel Rivi\`ere]{Gabriel Rivi\`ere}

\address{Universidad Polit\'ecnica de Madrid. ETSI Navales. Avda. de la Memoria, 4. 28040 Madrid, Spain}
\email{Fabricio.Macia@upm.es}
\thanks{FM takes part into the visiting faculty program of ICMAT and is
partially supported by grants ERC Starting Grant 277778 and
MTM2013-41780-P (MEC)}

\address{Laboratoire Paul Painlev\'e (U.M.R. CNRS 8524), U.F.R. de Math\'ematiques, Universit\'e Lille 1, 59655 Villeneuve d'Ascq Cedex, France}
\email{gabriel.riviere@math.univ-lille1.fr}
\thanks{GR is partially supported by the Agence Nationale de la Recherche through the
Labex CEMPI (ANR-11-LABX-0007-01) and the ANR project GeRaSic (ANR-13-BS01-
0007-01)}

\begin{abstract} 
We study the regularity of stationary and time-dependent solutions to strong perturbations of the free Schr\"odinger equation on two-dimensional flat tori. This is achieved by performing a second microlocalization 
related to the size of the perturbation and by analysing concentration and nonconcentration properties at this new scale. In particular, we show that sufficiently accurate quasimodes can only concentrate on the set of critical points of the average of the potential along geodesics. 
\end{abstract}

\maketitle

\section{Introduction}

The high-frequency analysis of eigenfunctions of elliptic operators on a compact Riemannian manifold has been the subject of intensive study in the past fifty years. To this day, many questions remain open, even in the simplest cases. Here we focus on eigenfunctions of Schrödinger operators on $\mathbb{T}^d:=\IR^d/\IZ^d$, the standard torus endowed with its canonical metric. 

Eigenfunctions of a Schrödinger operator on $\IT^d$ are precisely the solutions to the equation:
\begin{equation}\label{e:laplace-beltrami}
-\Delta u_{\lambda}(x)+V(x)u_{\lambda}(x)=\lambda^2 u_{\lambda}(x),\quad x\in\IT^d,\quad \quad \|u_{\lambda}\|_{L^2(\IT^d)}=1,
\end{equation}
where the potential $V$ is real-valued and essentially bounded. 
In the free case $V=0$, a straightforward computation shows that eigenfunctions of eigenvalue $\lambda^2$ are linear combinations of complex exponentials $e^{2i\pi k. x}$ with frequencies $k\in\IZ^d$ lying on a circle of radius $\lambda>0$ centered at the origin. However, extracting from this exact representation formula an asymptotic description of eigenfuctions in the high frequency limit $\lambda\rightarrow +\infty$ is a hard problem, due to the fact that multiplicities of large eigenvalues can also be very big. Instead, one can try to describe particular features of high-frequency eigenfunctions, such as formation of (asymtotic) singularities. 

A natural way to quantify these singularities is through the scale of $L^p$ spaces. This has been a classical topic in harmonic analysis, that originates with the seminal result of Zygmund~\cite{Zy74} showing that, for $d=2$ and in the free case, there exists some universal constant $C$ such that any solution $u_{\lambda}$ 
of~\eqref{e:laplace-beltrami} verifies $\|u_{\lambda}\|_{L^4(\IT^2)}\leq C$. Later on, Bourgain conjectured in~\cite{Bo93} that, again for the free case and when $d\geq 3$, one must have $\|u_{\lambda}\|_{L^{\frac{2d}{d-2}}(\IT^d)}\leq C_{\delta}\lambda^{\delta}$ for every $\delta>0$. We refer the reader 
to~\cite{Bo13, BoDe15} for recent progress towards this conjecture. Note that the problem of showing the existence of an index $p>2$ such that $\|u_{\lambda}\|_{L^p(\IT^d)}$ is uniformly bounded remains open for $d\geq 3$.

There are alternative ways to describe the asymptotic 
structure of the solutions of~\eqref{e:laplace-beltrami}. For instance, notice that a direct corollary of Zygmund's result is that, in the free case, any accumulation point of the sequence of probability measures,
$$\nu_{\lambda}(dx)=|u_{\lambda}(x)|^2dx,$$
is a probability measure which is absolutely continuous with respect to the Lebesgue measure on $\IT^2$ (it has in fact an $L^2$ density). This result was refined by Jakobson who 
showed that the density has to be a trigonometric polynomial whose frequencies enjoy certain geometric constraints~\cite{Ja97}. It is natural to try to understand what happens when $d\geq 3$, where no analogue to Zygmund's result is known to hold, or when the Laplacian is perturbed by a lower order term, such as a potential. Note that the problem of identifying accumulation points of sequences of moduli squares of eigenfunctions has a long history and it is connected to fundamental questions in quantum mechanics.

In dimension 
$d\geq 3$ and for $V=0$, Bourgain proved in \cite{Ja97} that any accumulation point has to be absolutely continuous even if we do not know \textit{a priori} that the $L^p$ 
norms of eigenfunctions are uniformly bounded for small $p>2$. In the same reference, Jakobson obtained partial results on the structure of the densities of accumulation points. These results are based on harmonic analysis techniques and arguments on the geometry of lattice points. Absolute continuity of accumulation points also holds in the case of a non-zero potential $V\in L^\infty(\IT^d)$, as was proved by Anantharaman and the first author~\cite{AM10}. The proof of that result is based on methods from semiclassical analysis for the time dependent Schrödinger equation that were introduced for the particular case $d=2$ in~\cite{Ma10}. In fact, the results in reference \cite{AM10} apply to the more general problem:
\begin{equation}\label{e:quasimode}
 \hat{P}_{\eps}(\hbar)u_{\hbar}=\frac{1}{2}u_{\hbar}+o(\hbar\eps_{\hbar}),\ \ \|u_{\hbar}\|_{L^2(\IT^d)}=1,
\end{equation}
where $\hbar\rightarrow 0^+$ is some semiclassical parameter, and where
\begin{equation}\label{e:schrodinger-op}
 \hat{P}_{\eps}(\hbar):=-\frac{\hbar^2\Delta}{2}+\eps_{\hbar}^2V,
\end{equation}
and with $0\leq \eps_{\hbar}\leq\hbar$ for $\hbar$ small enough.\footnote{Note that, when $\hbar=\eps_{\hbar}=\lambda^{-1}$, equation \eqref{e:quasimode} is essentially equation \eqref{e:laplace-beltrami}.}  One of the main ingredients used in this approach are the two-microlocal 
techniques developed in~\cite{Ni96, Mi96, Fe00, FeGe02} in a different context. The results in \cite{AM10} were further extended to treat the case of more general completely integrable 
systems in~\cite{AFM12}. Note that studying the regularity of the solutions to~\eqref{e:quasimode} is also related to problems arising in control theory as was 
shown by Burq and Zworski~\cite{BuZw04}. We refer the reader to~\cite{AL14, AM10, BoBuZw13, BuZw04, BuZw12,  MaciaDispersion} for perspectives 
from the point of view of control theory. 

A different but related approach consists in studying the wavefront set $WF_{\hbar}(u_{\hbar})$ of solutions to \eqref{e:quasimode}. This was done in series of works by Wunsch~\cite{Wu08, Wu12} and Vasy--Wunsch~\cite{VaWu09} dealing with completely integrable 
systems in dimension $d=2$. In these articles, the authors investigate the properties of the semiclassical wavefront set $WF_{\hbar}(u_{\hbar})$ of solutions 
to~\eqref{e:quasimode} when $0\leq\eps_{\hbar}\leq\hbar^{1+\delta}$ with $\delta>0$. By proving some propagation of second microlocal wavefront sets, they showed that $WF_{\hbar}(u_{\hbar})$ can 
not be reduced to a single geodesic and has to fill a Lagrangian torus -- see for instance~\cite[Th.~B]{Wu08} or~\cite[Th.~3]{Wu12}. Note that, 
as in~\cite{AFM12}, the results of Vasy and Wunsch hold for general classes of nondegenerate completely integrable systems. Under the assumption 
that $\hbar^{1-\delta}\ll\eps_{\hbar}\ll 1$, Wunsch also exhibited examples of quasimodes of 
order $\mathcal{O}(\hbar^{\infty})$ for the operator $\hat{P}_{\eps}(\hbar)$ which concentrate 
on closed geodesics; this result was reported in~\cite[Sect.~5.3]{AFM12}. This shows that $\eps_{\hbar}=\hbar$ is the critical size for which one can expect to have singular concentration phenomena for perturbations of the free 
semiclassical Schr\"odinger operator $-\frac{\hbar^2\Delta}{2}$. In particular, for stronger perturbation $\eps_{\hbar}\gg\hbar$, one cannot expect to have uniform bounds for $L^p$ norms 
even for small range of $p$. A notable feature of Wunsch's construction is that the singularity is located on critical points of the potential $V$ restricted 
to certain closed geodesics. In some sense, this type of singularities is similar to the ones that may occur in the case of Zoll manifolds~\cite{MaRi15, MaRi17}. Motivated by this 
observation, we will combine the ideas from~\cite{AM10, MaRi15} in order to derive some properties on the regularity of solutions 
to~\eqref{e:quasimode} when $\eps_{\hbar}\gg\hbar$. In particular, we will identify precisely the concentration phenomena that may occur and 
also show nonconcentration properties by propagation of second microlocal data.

For the sake of simplicity, we will focus on the case of the rational torus $\IT^2$ and assume that $V\in \ml{C}^{\infty}(\IT^2;\IR)$; but it is most likely that our analysis could be extended to more general completely 
integrable systems of dimension $2$ following the approach of~\cite{AFM12}. As the small perturbation regime $0\leq \eps_{\hbar}\leq\hbar$ was studied in great detail 
in all the above references, here we will focus on the strong perturbation regime and we shall assume all along the article that
\begin{equation}\label{e:size-perturbation}
 \lim_{\hbar\rightarrow 0^+}\eps_{\hbar}=0,\quad\text{and}\quad\lim_{\hbar\rightarrow 0^+}\hbar\eps_{\hbar}^{-1}=0.
\end{equation}
In order to state our results, we need some simple geometric preliminaries. Recall that the geodesics of $\IT^2$ are either closed or dense curves. For $\xi=(\xi_1,\xi_2)\in\IR^2-\{0\}$ and $x\in\IT^2$, the geodesic $s\mapsto x+s\xi$ is dense provided $\xi_1$ and $\xi_2$
are linearly independent over $\IQ$, otherwise it is periodic. We denote by $\Omega_1\subset \IR^2-\{0\}$ the set of $\xi$ that generate a periodic geodesic and by $\Omega_2$ its complementary in $ \IR^2-\{0\}$.  Consider the average of $V$ along geodesics:
$$\ml{I}(V)(x,\xi):=\lim_{T\to +\infty}\frac{1}{T}\int_0^T V(x+s\xi)ds.$$ 
Clearly, $\ml{I}(V)$ is a zero-homogeneous function with respect to $\xi$. Moreover, a classical result by Kronecker implies that
$$\ml{I}(V)(x,\xi)=\left\{\begin{array}{ll}
\frac{1}{L_\xi}\int_0^{L_\xi} V\left(x+s\frac{\xi}{\|\xi\|}\right)ds & \text{ if }\xi\in\Omega_1,\medskip \\
\int_{\IT^2}V(y)dy & \text{ if }\xi\in\Omega_2,
\end{array}\right. $$
where $L_\xi$ denotes the length of any geodesic with velocity $\xi$. In particular, $\ml{I}(V)(\cdot,\xi)\in \ml{C}^{\infty}(\IT^d;\IR)$ for any $\xi\in\IR^2-\{0\}$.

Then, we define the set of critical geodesics:
\begin{equation}\label{e:crit-geod}
 \ml{C}(V):=\left\{x_0\in \IT^2:\ \exists\ \xi\in\Omega_1\ \text{s.t.}\ \partial_x\ml{I}(V)(x_0,\xi)=0\right\}.
\end{equation}
Note that $\ml{C}(V)$ is a union of closed geodesics of $\IT^2$. For every closed geodesic $\gamma$ of $\IT^2$, 
we denote by $\delta_{\gamma}$ the normalized Lebesgue measure along this closed geodesic. Then, we define $\ml{N}(V)$ 
as the convex closure of the set of probability measures $\delta_{\gamma}$ where $\gamma\subset\ml{C}(V)$. With these conventions in mind, we can state our main result:
\begin{theo}\label{t:maintheo-quasimode} Suppose that $d=2$ and that~\eqref{e:size-perturbation} holds. Let $(u_{\hbar})_{\hbar\rightarrow 0^+}$ be a sequence satisfying~\eqref{e:quasimode}. Then, 
for any accumulation point $\nu$ of the sequence of 
probability measures
$$\nu_{\hbar}(dx):=|u_{\hbar}(x)|^2dx,$$
and for any closed geodesic $\gamma$, one has
$$\nu(\gamma)\neq 0\Longrightarrow\gamma\subset\ml{C}(V).$$
 Moreover, $\nu$ can be decomposed as
 $$\nu= fdx+\nu_{\operatorname{sing}},$$
 where $f\in L^1(\IT^2)$ and where $\nu_{\operatorname{sing}}\in\ml{N}(V)$.
\end{theo}
Recall from the propagation properties of semiclassical measures~\cite{Ge91, Zw12} that any such $\nu$ must a priori be a convex combination of the Lebesgue measure and of the 
measures $\delta_{\gamma}$, where $\gamma$ runs over the set of all closed geodesics. This Theorem shows that singular concentration along 
closed geodesics can only occur along certain closed orbits associated with critical points of the averages 
of $V$ along closed geodesics. This result is sharp in the sense that Wunsch's construction in~\cite{AFM12} shows that one can find quasimodes such that $\nu(\gamma)=1$ for a given closed 
geodesic. Despite these unavoidable concentration phenomena, Theorem~\ref{t:maintheo-quasimode} also shows that the accumulation points enjoy 
certain regularity properties. This extra regularity will come out from our analysis by making a second microlocalization of size $\eps_{\hbar}$ along 
rational directions, and, it will be induced by certain Lagrangian tori associated to our problem. Note that these two aspects are close to the situation of Zoll manifolds 
treated in~\cite{MaRi15, MaRi17}. The main difference is that there exist infinitely many directions where the flow is periodic with periods tending 
to $+\infty$. We would like to treat these tori of periodic orbits as in this reference, and this can be achieved via rescaling the variables along these rational directions -- see paragraph~\ref{ss:zoll} for more details. Finally, 
as we shall see it in section~\ref{s:timedependent}, our analysis holds in the more 
general context of the time dependent Schr\"odinger equation.

\subsection*{Organization of the article} 

In section~\ref{s:2-micro-statement}, we introduce the $2$-microlocal framework of our analysis and formulate our main results using this terminology. In 
section~\ref{s:timedependent}, we show how to apply the results of section~\ref{s:2-micro-statement} in order to study the semiclassical measures of the time 
dependent Schr\"odinger equation associated with $\hat{P}_{\eps}(\hbar)$ and in order to derive Theorem~\ref{t:maintheo-quasimode}. The proofs of the $2$-microlocal 
statements is given in section~\ref{s:2microlocal}. Finally, the article contains two appendices. Appendix~\ref{a:biinvariant} contains the proof of a geometric 
result which already appeared in~\cite{MaRi15} and which we adapt to the context of $\IT^2$. In Appendix~\ref{a:sc-an}, we collect a few tools from semiclassical 
analysis.

In the following (except in appendix~\ref{a:sc-an}), we will always suppose that $d=2$ and that~\eqref{e:size-perturbation} holds even if part of the results 
holds in greater generality.

\section{Invariance and propagation of $2$-microlocal distributions}\label{s:2-micro-statement}

As was already mentionned, Theorem~\ref{t:maintheo-quasimode} is a consequence of our analysis of the time dependent semiclassical Schr\"odinger equation:
\begin{equation}\label{e:perturbed-schr}
i\hbar\partial_t v_{\hbar}=\hat{P}_{\eps}(\hbar)v_{\hbar},\quad v_{\hbar}|_{t=0}=u_{\hbar}\in L^2(\IT^2),\quad\|u_{\hbar}\|_{L^2}=1.
\end{equation}
More specifically, our main results describe the $2$-microlocal structure of these solutions along covectors in $\Omega_1$. In other words, we will show how solutions of~\eqref{e:perturbed-schr} can concentrate along rational 
covectors. Let us now be more precise. 

\subsection{Set-up and conventions}

First of all, we shall focus for the sake of simplicity on sequences of initial data oscillating at 
the frequency $\hbar^{-1}$. Thus, we will always assume that the following properties hold:
\begin{equation}\label{e:hosc}
\limsup_{\hbar\rightarrow 0}\left\Vert \mathbf{1}_{\left[  R,\infty\right)
}\left(  -\hbar^{2}\Delta\right)  u_{\hbar}\right\Vert _{L^{2}\left(  M\right)
}\longrightarrow0,\quad\text{as }R\longrightarrow\infty,
\end{equation}
and%
\begin{equation}\label{e:shosc}
\limsup_{\hbar\rightarrow 0}\left\Vert \mathbf{1}_{\left[  0,\delta\right]
}\left(  -\hbar^{2}\Delta\right)  u_{\hbar}\right\Vert _{L^{2}\left(  M\right)
}\longrightarrow0,\quad\text{as }\delta\longrightarrow0^{+}.
\end{equation}
For every primitive rank $1$ lattice\footnote{This just means that $\text{dim}\la\Lambda\ra=1$ and that 
$\la \Lambda\ra\cap \IZ^2=\Lambda$, where $\la\Lambda\ra$ is the linear subspace of $\IR^2$ spanned by $\Lambda$.} $\Lambda$ of $\IZ^2$, we set $\mathfrak{e}_{\Lambda}$ to be an element in $\Lambda$ such that $\IZ \mathfrak{e}_{\Lambda}=\Lambda$, 
and $\mathfrak{e}_{\Lambda}^{\perp}$ to be the vector of same length which is directly orthogonal to $\mathfrak{e}_{\Lambda}$. We define 
$$L_{\Lambda}:=\|\mathfrak{e}_{\Lambda}\|.$$
We define two Hamiltonian maps associated to $\Lambda$ as follows:
$$H_{\Lambda}(x,\xi):=\frac{1}{L_{\Lambda}}\la\xi,\mathfrak{e}_{\Lambda}\ra\ \text{and}\ H_{\Lambda}^{\perp}(x,\xi):=\frac{1}{L_{\Lambda}}\la\xi,\mathfrak{e}_{\Lambda}^{\perp}\ra.$$
Note that $(H_{\Lambda},H_{\Lambda}^{\perp})$ defines a (nondegenerate) completely integrable system and that 
$$\|\xi\|^2=H_{\Lambda}(x,\xi)^2+H_{\Lambda}^{\perp}(x,\xi)^2.$$

\subsection{Two-microlocal distributions}

We aim at studying the concentration of solutions to \eqref{e:perturbed-schr} over $\IT^2\times\Lambda^\bot$ where $\Lambda\subset\IZ^2$ is a primitive rank $1$ sublattice and where $\Lambda^{\perp}$ denotes the set of covectors $\xi$ such that $\la\xi,\mathfrak{e}_{\Lambda}\ra=0$. For that purpose, we define the following two-microlocal Wigner distribution:
$$w_{\Lambda,\hbar}(t):a\in\ml{C}^{\infty}_c(T^*\IT^2\times\widehat{\IR})\longmapsto \left\la v_{\hbar}(t),
\Oph^w\left(a\left(x,\xi,\frac{H_{\Lambda}(x,\xi)}{\eps_{\hbar}}\right)\right)v_{\hbar}(t)\right\ra.$$
Above, $\widehat{\IR}$ is the compactified space $\IR\cup\{\pm\infty\}$, $v_{\hbar}(t)$ is the solution of~\eqref{e:perturbed-schr} 
at time $t$ and $\Oph^w(a)$ is a $\hbar$-pseudodifferential operator -- see Appendix~\ref{a:sc-an}. 
\begin{rema}\label{r:semiclas}
Recall from~\eqref{e:change-variable} in Appendix~\ref{a:sc-an} that the following useful relation holds:
$$\Oph^w\left(a\left(x,\xi,\frac{H_{\Lambda}(x,\xi)}{\eps_{\hbar}}\right)\right)=\Op_{\hbar\eps_{\hbar}^{-1}}^w
\left(a\left(x,\eps_{\hbar}\xi,H_{\Lambda}(x,\xi)\right)\right),$$
and that we made the assumption that $\hbar\eps_{\hbar}^{-1}\rightarrow 0$. Therefore, the operators involved in the definition of $w_{\Lambda,\hbar}$ are semiclassical pseudodifferential operators.
\end{rema}

\begin{rema}
The distributions $w_{\Lambda,\hbar}$ were introduced in \cite{Ma10, AM10} for the critical case $\eps_\hbar = \hbar$. As we will see, its limiting objects are of a very different nature in the present case.
\end{rema}
 
Fix now a sequence of time scales $(\tau_{\hbar})_{\hbar\rightarrow 0^+}$ such that
$$\lim_{\hbar\rightarrow 0^+}\tau_{\hbar}=+\infty.$$
As we shall explain it in paragraph~\ref{ss:extract}, we can extract a subsequence $\hbar_n\rightarrow 0^+$ such that, for any 
$a\in\ml{C}^{\infty}_c(T^*\IT^2\times\widehat{\IR})$ and for any $\theta\in L^1(\IR)$,
$$\lim_{n\rightarrow +\infty}\int_{\IR}\theta(t)\la w_{\Lambda,\hbar_n}(t\tau_{\hbar_n}),a\ra dt=\int_{\IR}\theta(t)
\left(\int_{T^*\IT^2\times\widehat{\IR}}a(x,\xi,\eta)\mu_{\Lambda}(t,dx,d\xi,d\eta)\right)dt,$$
where, for a.e. $t$ in $\IR$, $\mu_{\Lambda}(t)$ is an element of $\ml{B}'$ for some Banach space $\ml{B}$ that we will define 
in paragraph~\ref{ss:extract}. We denote by $\ml{M}_{\Lambda}(\tau,\eps)$ the set of accumulation points obtained in this manner 
for initial data varying among subsequences verifying~\eqref{e:hosc} and~\eqref{e:shosc}. The main new result of this article describes 
some invariance and propagation properties of these quantities depending on the relative size of $\tau_{\hbar}$ and 
$\eps_{\hbar}$.

Before stating our main results, we will show that any element $\mu_{\Lambda}(t)$ inside $\ml{M}(\tau,\eps)$ is a measure that is concentrated on $\mathring{T}^*\IT^2\times\widehat{\IR} $, where
$$\mathring{T}^*\IT^2:=\left\{(x,\xi)\in T^*\IT^2: \xi\neq 0\right\}.$$ 
We will split them
in two components:
\begin{equation}\label{e:split-distrib}
 \mu_{\Lambda}(t)=\tilde{\mu}_{\Lambda}(t)+\tilde{\mu}^{\Lambda}(t).
\end{equation}
with $\tilde{\mu}_{\Lambda}(t)$ corresponding to the restriction to the ``finite'' part $T^*\IT^2\times\IR$ and $\tilde{\mu}^{\Lambda}(t)$ to the part at 
infinity $T^*\IT^2\times\{\pm\infty\}$. Hence, $\tilde{\mu}_{\Lambda}(t)$ describes 
in some sense the way the solutions of~\eqref{e:perturbed-schr} concentrate in an $\eps_{\hbar}$-neighborhood of the rational direction 
$\Lambda^\bot$. Let us start by giving some simple properties of these functionals which are analoguous to the ones satisfied by 
time dependent semiclassical measures~\cite{Ma09}. 
\begin{prop}\label{p:first-properties} Let $\mu_{\Lambda}(t)$ be an element of $\ml{M}_{\Lambda}(\tau,\eps)$. Then, for a.e. $t$ in $\IR$, $\mu_{\Lambda}(t)$ is a positive finite Radon measure concentrated on $\mathring{T}^*\IT^2\times\widehat{\IR} $. Let $$\tilde{\mu}_\Lambda(t):=\mu_\Lambda (t)\rceil_{T^*\IT^2\times \IR},\quad \tilde{\mu}^\Lambda(t):=\mu_\Lambda (t)\rceil_{T^*\IT^2\times \{\pm \infty\}},$$
so that \eqref{e:split-distrib} holds. Then
\begin{enumerate}
 \item $\tilde{\mu}_{\Lambda}(t)$ is a (finite) positive measure on $T^*\IT^2\times\IR$ whose support is contained in $\IT^2\times(\Lambda^{\perp}-\{0\})\times\IR$;
 \item for every $a$ in $\ml{C}_c^{\infty}(T^*\IT^2\times\widehat{\IR})$,
 $$\la\tilde{\mu}_{\Lambda}(t),\xi.\partial_xa\ra=\la\tilde{\mu}^{\Lambda}(t),\xi.\partial_xa\ra=0.$$
\end{enumerate}
\end{prop}

These properties follow from standard arguments which need to be slightly adapted in order to fit into the $2$-microlocal 
set-up -- see Section~\ref{s:2microlocal} for details.

\subsection{Main results} 

Consider the Hamiltonian flow 
$\varphi_{H_{\Lambda}^{\perp}}$ associated with $H_{\Lambda}^{\perp}$. Note that, for a continuous function $b$ on $T^*\IT^2\times\widehat{\IR}$, we 
can define the average along this $L_\Lambda$-periodic flow as
$$\ml{I}_{\Lambda}(b)(x,\xi,\eta):=\frac{1}{ L_{\Lambda}}\int_0^{ L_{\Lambda}}b\left(\varphi_{H_{\Lambda}^{\perp}}^s(x,\xi),\eta\right)ds.$$
A direct computation gives
$$\ml{I}_{\Lambda}(b)(x,\xi,\eta)=\frac{1}{L_\Lambda}\int_0^{L_\Lambda}b\left(x+s\frac{\mathfrak{e}_{\Lambda}^{\perp}}{L_\Lambda},\xi,\eta\right)ds=\sum_{k\in \Lambda}\hat{b}_k(\xi,\eta)e^{2i\pi k.x},$$
provided $b$ has the Fourier expansion $b(x,\xi,\eta)=\sum_{k\in\IZ^2}\hat{b}_k(\xi,\eta)e^{2i\pi k.x}$. Moreover, if $\ml{I}(b)$ denotes the average of $b$ along the geodesic flow 
$$\varphi^s(x,\xi)=(x+s\xi,\xi)$$ 
on $T^*\IT^2$ then the following holds:
\begin{equation}\label{e:eqflow}
\ml{I}(b)(x,\xi,\eta)=\ml{I}_{\Lambda}(b)(x,\xi,\eta),\quad \text{ provided that }\xi\in\Lambda^\bot-\{0\}.
\end{equation}

\begin{rema}\label{r:invmul}
Part (2) of Proposition \eqref{p:first-properties} implies that $\mu_\Lambda(t)$ is invariant under the geodesic flow $\varphi^s$. For $b$ in $\ml{C}^{\infty}_c(T^*\IT^2\times\IR)$, this observation combined with part (1) in Proposition~\ref{p:first-properties} and identity \eqref{e:eqflow} implies that, for a.e. 
 $t$ in $\IR$,
 $$\la \mu_{\Lambda}(t),b\ra=\la\mu_{\Lambda}(t),\ml{I}_{\Lambda}(b)\ra.$$
 We shall use this property several times in our proof of Theorem~\ref{t:propag-inv} below.
\end{rema}
In the case where $b$ only depends on $x$, as is the case with $b=V$, it is easy to check that $\ml{I}_{\Lambda}(V)$ does not depend on $\xi$ and therefore we can identify it to an element in $ \ml{C}^{\infty}(\IT^d;\IR)$.

We need to define an auxiliary Hamiltonian function on $\IT^2\times\Lambda^{\perp}\times\IR$
\begin{equation}\label{e:comm-ham} p_{\Lambda}(x,\sigma\mathfrak{e}_\Lambda^\bot/L_\Lambda,\eta):=\frac{\eta^2}{2}+\ml{I}_{\Lambda}(V)(x).\end{equation}
Denote by $\varphi_{p_{\Lambda}}^t$ the flow of the vector field on $\IT^2\times\Lambda^{\perp}\times\IR$: 
$$\eta\frac{\mathfrak{e}_{\Lambda}}{L_{\Lambda}}.\partial_x-\frac{\mathfrak{e}_{\Lambda}}{L_{\Lambda}}.\partial_x\ml{I}_{\Lambda}(V)\partial_{\eta}.$$
This is the Hamiltonian vector field associated to $p_{\Lambda}$ with respect to the symplectic form
obtained by taking the push-forward of the canonical symplectic form on $T^*\IT^2$ via the diffeomorphism
\begin{equation}\label{e:diffeo}
T^*\IT^2\ni (x,\xi) \longmapsto (x,H_\Lambda^\bot(x,\xi)\mathfrak{e}_\Lambda^\bot/L_\Lambda, H_\Lambda(x,\xi))\in \IT^2\times\Lambda^{\perp}\times\IR.
\end{equation}
The flow $\varphi_{p_{\Lambda}}^t$ commutes with $\varphi_{H_{\Lambda}^{\perp}}^s$ when acting on $\IT^2\times\Lambda^{\perp}\times\IR$.

We are now ready to state the main results of this article. 
The first one concerns the ``compact'' part of these two-microlocal distributions:
\begin{theo}[Invariance and propagation near $\Lambda$]\label{t:propag-inv} Let $\Lambda$ be a primitive rank $1$ sublattice and let 
$\mu_{\Lambda}$ be an element of $\ml{M}_{\Lambda}(\tau,\eps)$ obtained as the limit of $(w_{\Lambda,\hbar})$. Denote by $\tilde{\mu}_{\Lambda}^0$ the limit of $(w_{\Lambda,\hbar}(0))$. The following results hold:
\begin{enumerate}
  \item If $\tau_{\hbar}\eps_{\hbar}\rightarrow 0$ as $\hbar\rightarrow 0^+$, then $t\mapsto \tilde{\mu}_{\Lambda}(t)$ is continuous, 
  and one has, for every $a$ in $\ml{C}^{0}_c(\IT^2\times\Lambda^{\perp}\times\IR)$,
  $$\tilde{\mu}_{\Lambda}(t)(a)=\tilde{\mu}_{\Lambda}^0(\ml{I}_{\Lambda}(a)).$$
  \item If $\tau_{\hbar}\eps_{\hbar}=1$ as $\hbar\rightarrow 0^+$, then $t\mapsto \tilde{\mu}_{\Lambda}(t)$ is continuous, and one has, 
  for every $a$ in $\ml{C}^{0}_c(\IT^2\times\Lambda^{\perp}\times\IR)$,
  $$\tilde{\mu}_{\Lambda}(t)(a)=\tilde{\mu}_{\Lambda}^0(\ml{I}_{\Lambda}(a)\circ\varphi_{p_{\Lambda}}^t).$$
  \item If $\tau_{\hbar}\eps_{\hbar}\rightarrow +\infty$ as $\hbar\rightarrow 0^+$, then one has, for a.e. $t$ in $\IR$ and, 
  for every every $a$ in $\ml{C}^{0}_c(\IT^2\times\Lambda^{\perp}\times\IR)$,
  $$\forall s\in\IR,\quad\tilde{\mu}_{\Lambda}(t)\left(a\right)=\tilde{\mu}_{\Lambda}(t)\left(a\circ\varphi_{p_{\Lambda}}^s\right).$$
\end{enumerate}
\end{theo}
Equivalently, this Theorem says that, besides invariance by the geodesic flow, the solutions of~\eqref{e:perturbed-schr} 
satisfy some extra invariance properties in a shrinking neighborhood of the rational direction at least for times 
$\tau_{\hbar}\gg\frac{1}{\eps_{\hbar}}$. For shorter times, the concentration in this shrinking neighborhood is completely 
determined by the initial data.

For the part at infinity, we have the following 
regularity property:
\begin{theo}[Regularity at infinity]\label{t:maintheo-infty} Let $\Lambda$ be a primitive rank $1$ sublattice and let 
$\mu_{\Lambda}(t)$ be an element of $\ml{M}_{\Lambda}(\tau,\eps)$. Suppose that
$\tau_{\hbar}\eps_{\hbar}\geq 1$;
then, one has, for every $k\in \Lambda-\{0\}$, for every $a$ in $\ml{C}^{\infty}_c(\IR^2\times\widehat{\IR})$ and for a.e. $t$ in $\IR$,
$$\la\tilde{\mu}^{\Lambda}(t),a(\xi,\eta)e^{-2i\pi k.x}\ra=0.$$
In particular, the measure $\tilde{\mu}^{\Lambda}(t)\rceil_{\IT^2\times \Lambda^\bot\times \widehat{\IR}}$ is constant in $x$.
\end{theo}
In other words, the part at infinity has no (nonzero) Fourier coefficients in the $\Lambda$-direction for large enough scales of times. 
The last assertion of the Theorem follows from the invariance\footnote{Recall also that $\mu_{\Lambda}(t)$ is supported on $\mathring{T}^*\IT^2\times\widehat{R}$.} 
of $\tilde{\mu}^{\Lambda}(t)$ under the geodesic flow, 
which implies that for every $a\in \ml{C}^0_c(T^*\IT^2\times\widehat{\IR}) $:
$$\la\tilde{\mu}^{\Lambda}(t)\rceil_{\IT^2\times \Lambda^\bot\times \widehat{\IR}},a\ra=
\la\tilde{\mu}^{\Lambda}(t)\rceil_{\IT^2\times \Lambda^\bot\times \widehat{\IR}},\ml{I}_\Lambda (a)\ra= 
\la\tilde{\mu}^{\Lambda}(t)\rceil_{\IT^2\times \Lambda^\bot\times \widehat{\IR}},\int_{\IT^2}a(y,\cdot)dy\ra,$$
since $\ml{I}_\Lambda (a) $ has only Fourier coefficients in the $\Lambda$-direction.

It is interesting to compare these result with the corresponding ones in \cite{AM10}, particularly with Corollary 25 in that reference. The propagation law in the critical case $\eps_\hbar=\hbar$ involves a quantum flow rather than a classical one .
\subsection{Comparison with Zoll manifolds}\label{ss:zoll}

Theorem~\ref{t:propag-inv} shares a lot of similarities with our main result on semiclassical measures 
for perturbations of Zoll Laplacians in~\cite[Sect.~2.2]{MaRi15}. In that case, we were considering the semiclassical operator
$$-\frac{\hbar^2\Delta_g}{2}+\eps_{\hbar}^2V,$$
 where $\Delta_g$ is the Laplace Beltrami operator associated to a certain Zoll metric 
 (say the standard metric on the canonical sphere). In the present article, we are analyzing the semiclassical measures associated to 
 the same Schr\"odinger operator $\hat{P}_{\eps}(\hbar)$. Studying the ``compact'' part of elements inside $\ml{M}_{\Lambda}(\tau,\eps)$ 
 is equivalent to understanding the solutions of~\eqref{e:perturbed-schr} near submanifolds
$$\IT^2\times\Lambda^{\perp}:=\{(x,\xi)\in T^*\IT^2:H_{\Lambda}(x,\xi)=0\},$$
where the geodesic flow is periodic as in the Zoll case. In order to make the comparison more clear and to justify the rescaling of order $\eps_{\hbar}$, 
we can rewrite our operator in a form which is very close to what we did in the Zoll framework, i.e.
$$\hat{P}_{\eps}(\hbar)=\frac{1}{2}\Oph^w(H_{\Lambda}^{\perp})^2+\eps_{\hbar}^2
\Op_{\hbar}^w\left(\frac{1}{2}\left(\frac{H_{\Lambda}}{\eps_{\hbar}}\right)^2+V\right).$$ 
Thus, as in the Zoll case, we perturb in some sense a semiclassical operator $\Oph^w(H_{\Lambda}^{\perp})^2$ asssociated to a ``periodic'' Hamiltonian 
vector flow and we obtain limit quantities which are invariant by the periodic flow and the Hamiltonian perturbation. 

The main difference with the Zoll setting is that the perturbation depends on rescaled variables  
$$\left(x,H_{\Lambda}^{\perp}(x,\xi),\frac{H_{\Lambda}(x,\xi)}{\eps_{\hbar}}\right)\in \IT^2\times\IR^2\simeq T^*\IT^2.$$
For that reason, it is natural to test our Wigner distributions against symbols depending on these rescaled variables. Another notable 
difference with~\cite{MaRi15} is that, in the Zoll case, the critical time scale is of order $\eps_{\hbar}^{-2}$ 
while here, due to the use of rescaled variables, it is much shorter, i.e. of order $\eps_{\hbar}^{-1}$. Finally, in the Zoll case, 
a natural question was to discuss the case where the Radon transform of the perturbation identically vanishes~\cite{MaRi17}. 
Here, we emphasize that the $H_{\Lambda}^{\perp}$-average of the perturbation, namely 
$\frac{1}{2}\left(\frac{H_{\Lambda}}{\eps_{\hbar}}\right)^2+\ml{I}_{\Lambda}(V)$ 
\emph{cannot be equal to a constant} for this choice of $2$-microlocal rescaling.

\section{From Theorems~\ref{t:propag-inv} and~\ref{t:maintheo-infty} to Theorem~\ref{t:maintheo-quasimode}}\label{s:timedependent}

Before proving our results on $2$-microlocal regularity, we show how to derive Theorem~\ref{t:maintheo-quasimode} from these results. In fact, we will prove something slightly stronger related to the time dependent semiclassical measures associated 
with the semiclassical Schr\"odinger equation~\eqref{e:perturbed-schr}. 

\subsection{Time-dependent semiclassical measures}

For a given $t$ in $\IR$, we denote the Wigner distribution at time $t$ by
\begin{equation}\label{e:wigner}w_{\hbar}(t)(a):=\left\la v_{\hbar}(t),\Oph^w(a)v_{\hbar}(t)\right\ra,
\end{equation}
where $\Oph^w(a)$ is a $\hbar$-pseudodifferential operator with principal symbol $a\in \ml{C}^{\infty}_c(T^*\IT^2)$ -- see Appendix~\ref{a:sc-an}. Again 
$v_{\hbar}(t)$ denotes the solution at time $t$ of~\eqref{e:perturbed-schr} with initial conditions satisfying the oscillating 
assumptions~\eqref{e:hosc} and~\eqref{e:shosc}. Observe that this is just a particular case of the two-microlocal distributions we have already 
introduced. This quantity represents the distribution of the $L^2$-mass of the solution to~\eqref{e:perturbed-schr} in the phase space $T^*\IT^2$. According 
to~\cite{Ma09}, we can extract a subsequence $\hbar_n\rightarrow 0^+$ as $n\rightarrow +\infty$ such that, for every $a$ in $\ml{C}^{\infty}_c( T^*\IT^2)$ 
and for every $\theta$ in $L^1(\IR)$,
$$\lim_{\hbar_n\rightarrow 0^+}\int_{\IR\times T^*\IT^2}\theta(t)a(x,\xi) w_{\hbar_n}(t\tau_{\hbar_n},dx,d\xi)dt=\int_{\IR\times T^*\IT^2}\theta(t)a(x,\xi)\mu(t,dx,d\xi)dt,$$
where, for a.e. $t$ in $\IR$, $\mu(t)$ is a finite positive Radon measure
on $T^*\IT^2$. Recall also that, for a.e. $t\in\IR$, $\mu(t)$ is in fact a \emph{probability measure} which does not put any mass on the zero section, thanks to the frequency assumption~\eqref{e:shosc}. In other words, 
\begin{equation}\label{e:tmass}
\mu(t)(\mathring{T}^*\IT^2)=1, \text{ for a.e. } t\in\IR.
\end{equation}
Moreover, for a.e. $t$ in $\IR$, $\mu(t)$ is \emph{invariant by the geodesic flow} $\varphi^s$ on $T^*\IT^2$. 

For instance, $\mu(t)$ can be the normalized Lebesgue measure along a closed orbit of the geodesic flow. We will denote by $\ml{M}(\tau,\eps)$ the set of accumulation
points of the sequences $(\mu_{\hbar})$, where 
$\mu_{\hbar}(t,\cdot):=
w_{\hbar}(t\tau_{\hbar},\cdot)$, as the sequence of initial data $(u_{\hbar})$ 
varies among normalized sequences satisfying~\eqref{e:hosc} and~\eqref{e:shosc}. For every primitive rank $1$ sublattice one has (see Remark \ref{r:proj}),
\begin{equation}\label{e:pushforward-2microlocal}
\ml{M}(\tau,\eps)=\left\{\int_{\widehat{\IR}}\mu_{\Lambda}(x,\xi,d\eta):\ \mu_{\Lambda}\in\ml{M}_{\Lambda}(\tau,\eps)\right\}.
 \end{equation}
Similarly, one can define $\ml{N}(\tau,\eps)$ to be the set of accumulation points of the sequences $(n_\hbar)$ of time-dependent probability measures on $\IT^2$, 
$n_\hbar(t,dx):=|v_{\hbar}(t\tau_{\hbar},x)|^2dx$, obtained letting the initial data vary among sequences satisfying~\eqref{e:hosc}, \eqref{e:shosc}. Using~\eqref{e:hosc}, one can verify that
\begin{equation}\label{e:pushforward}
\ml{N}(\tau,\eps)=\left\{\int_{\IR^2}\mu(x, d\xi):\ \mu\in\ml{M}(\tau,\eps)\right\}.
\end{equation}
In order to relate this to the quasimode case, we can remark that, given a sequence of quasimodes $(u_{\hbar})_{\hbar\rightarrow 0^+}$ 
satisfying~\eqref{e:quasimode}, we can always find a sequence of time scales $(\tau_{\hbar})$ such that
$$\lim_{\hbar\rightarrow 0}\tau_{\hbar}\eps_{\hbar}^{-1}=+\infty,$$
and, for every $t\in\IR$:
$$\lim_{\hbar\to 0}\| v_\hbar(\tau_\hbar t,\cdot) - e^{-i\tau_\hbar t /2\hbar} u_\hbar\|_{L^2(\IT^2)}=0,$$
where $v_\hbar$ denotes the solution to \eqref{e:perturbed-schr} with initial condition $u_\hbar$. This choice of $(\tau_\hbar)$ ensures that any accumulation point $\nu$ of the sequence of probability measures $(|u_{\hbar}|^2dx)$ belongs to $\ml{N}(\tau,\eps)$ (even though it is constant in $t$), since it is also an accumulation point of $(|v_{\hbar}(\tau_\hbar t,\cdot)|^2dx)$. In particular, Theorem~\ref{t:maintheo-quasimode} follows from the more general statement:
\begin{theo}\label{t:maintheo-schrodinger} Suppose that
$$\lim_{\hbar}\tau_{\hbar}\eps_{\hbar}^{-1}=+\infty.$$
Let $t\longmapsto\nu(t)$ be an element of $\ml{N}(\tau,\eps)$. Then, for any closed geodesic $\gamma$ not included inside $\ml{C}(V)$ and for a.e. $t$ in $\IR$, 
one has
$$\nu(t)(\gamma)= 0.$$
 Moreover, $\nu(t)$ can be decomposed as
 $$\nu(t)= f(t)dx+\nu_{\operatorname{sing}}(t),$$
 where, for a.e. $t$ in $\IR$, $f(t)\in L^1(\IT^2)$ and $\nu_{\operatorname{sing}}(t)\in \ml{N}(V)$.
\end{theo}

\subsection{Proof of Theorem~\ref{t:maintheo-schrodinger}} Let $t\longmapsto\mu(t)$ be an element of $\ml{M}(\tau,\eps)$. We start by splitting 
$\IR^2-\{0\}$ into $\varphi^s$-invariant subsets in the following manner. We introduce the set of rational covectors
$$\Omega_1=\bigsqcup_{\Lambda\ \text{rank 1 primitive}}\Lambda^{\perp}-\{0\},$$
and its complement $\Omega_2$ inside $\IR^2-\{0\}$. Observe that this is consistent with the conventions of the introduction. Because of \eqref{e:tmass},
we can decompose the measure as follows:
$$\mu(t)=\mu(t)\rceil_{\IT^2\times\Omega_2}+\sum_{\Lambda\ \text{rank 1 primitive}}\mu(t)\rceil_{\IT^2\times\Lambda^{\perp}-\{0\}}.$$
From the invariance by the geodesic flow, it can be verified that $\mu(t)\rceil_{\IT^2\times\Omega_2}$ is in fact independent of the 
$x$-variable. Hence, in order to prove Theorem~\ref{t:maintheo-schrodinger}, it remains to study the regularity of 
$\mu(t)\rceil_{\IT^2\times\Lambda^{\perp}-\{0\}}$ for every rank $1$ primitive sublattice $\Lambda$. This is where we will use our 
two-microlocal results. 
Thanks to~\eqref{e:pushforward-2microlocal} and to Proposition~\ref{p:first-properties}, we deduce
$$\mu(t)\rceil_{\IT^2\times\Lambda^{\perp}-\{0\}}=\mu(t)\rceil_{\IT^2\times\Lambda^{\perp}}=\int_{\IR}\tilde{\mu}_{\Lambda}(t,\cdot, d\eta)\rceil_{\IT^2\times\Lambda^{\perp}}
+\int_{\{\pm\infty\}}\tilde{\mu}^{\Lambda}(t,\cdot,d\eta)\rceil_{\IT^2\times\Lambda^{\perp}}.$$
According to Theorem~\ref{t:maintheo-infty}, the contribution from the part at infinity is independent of $x$. Hence, we are left 
with studying the regularity of the measures on $\IT^2$:
$$\int_{\Lambda^{\perp}\times\IR}\tilde{\mu}_{\Lambda}(t,\cdot,d\xi,d\eta).$$
The measure $\tilde{\mu}_{\Lambda}$
is invariant under the Hamiltonian flow $\varphi_{H_{\Lambda}^{\perp}}$ (see Remark \ref{r:invmul}) and, by  part~(3) of Theorem~\ref{t:propag-inv}, it is also invariant under the Hamiltonian flow $\varphi_{p_\Lambda}^t$, which commutes with $\varphi_{H_{\Lambda}^{\perp}}$. Using Appendix~\ref{a:biinvariant} 
which describes the regularity of biinvariant measures, we can conclude the proof of Theorem~\ref{t:maintheo-schrodinger}. More specifically, part 1 follows from 
Proposition~\ref{p:biinv} and part 2 from Corollary~\ref{c:struct}.

\subsection{Study of the critical time scale $\tau_{\hbar}=\eps_{\hbar}^{-1}$}

Let us now discuss what happens at the critical time scale
$$\tau_{\hbar}=\frac{1}{\eps_{\hbar}}.$$
In that case, it turns out that that the semiclassical measure can be completely determined from the initial data used to generate $\mu(t)$. 
More specifically, if we set $\tilde{\mu}_{\Lambda}^0$ to be the ``compact'' part of the two-microlocal distribution associated with the initial data 
and $\mu^0$ to be the semiclassical measure of the sequence of initial data, then $\mu(t)$ can be explicitly written in terms of these quantities. For 
that purpose, we shall start by recalling the following Lemma from~\cite[Prop.~29]{AM10}:
\begin{lemm}\label{l::initialdata} Suppose that
$$\lim_{\hbar\rightarrow 0^+}\tau_{\hbar}\eps_{\hbar}^2=0.$$
Let $\mu$ be an element in $\ml{M}(\tau,\eps)$ and let $\mu^0$ be the semiclassical
measure of the sequence of initial data used to generate $\mu$. Then, one has, for a.e. $t$ in $\IR$,
$$\int_{\IT^2}\mu(t,dy,\xi)=\int_{\IT^2}\mu^0(dy,\xi).$$
\end{lemm}
Arguing as before, if we fix $\mu(t)$ in $\ml{M}(\tau,\eps)$, then we can decompose it in three parts as follows
\begin{align*}
\mu(t)  =  \mu(t)\rceil_{\IT^2\times\Omega_2}+& \sum_{\Lambda\ \text{rank 1 primitive}}
\int_{\{\pm\infty\}}\tilde{\mu}^{\Lambda}(t,d\eta)\rceil_{\IT^2\times\Lambda^{\perp}}\\+&\sum_{\Lambda\ \text{rank 1 primitive}}\int_{\IR}\tilde{\mu}_{\Lambda}(t,\cdot,d\eta)\rceil_{\IT^2\times\Lambda^{\perp}}.
\end{align*}
Thanks to the invariance by the geodesic flow and to Theorem~\ref{t:maintheo-infty}, we can conclude one more time that the first two terms on the 
right-hand side of the equality are independent of $x$. Thanks to the second part of Theorem~\ref{t:propag-inv}, we can also write:
$$\tilde{\mu}_{\Lambda}(t)\rceil_{\IT^2\times\Lambda^{\perp}\times\IR}
=\left(\varphi_{p_{\Lambda}}^{t}\right)_*\left(\tilde{\mu}_{\Lambda}^0\rceil_{\IT^2\times\Lambda^{\perp}\times\IR}\right).$$
Hence, it is completely determined by the initial data. As the zero Fourier coefficient of $\mu(t)$ is itself equal to the zero Fourier 
coefficient of $\mu^0$ thanks to Lemma~\ref{l::initialdata}, we finally find that $\mu(t)$ can be expressed only in terms of the initial data.

\section{Proof of the $2$-microlocal statements}
\label{s:2microlocal}

From this point on, we fix a primitive lattice $\Lambda$ of $\IZ^2$ of rank $1$ and we will proceed to the proofs of the results on 
$2$-microlocal distributions. Namely, we will first recall how to extract converging subsequences from the sequences 
$(w_{\Lambda,\hbar})_{\hbar\rightarrow 0^+}$. Then, we will briefly recall how to adapt the proofs from~\cite{AM10} in order to prove 
Proposition~\ref{p:first-properties}. Finally, we will give the proofs of Theorems~\ref{t:propag-inv} and~\ref{t:maintheo-infty}.

\subsection{Extracting subsequences}\label{ss:extract}

Recall that, following~\cite{Ma10, AM10, AFM12}, we have introduced an auxiliary linear form whose invariance properties will be analyzed precisely. 
For every $a\in\ml{C}^{\infty}_c(T^*\IT^2\times\widehat{\IR})$, we set
$$\la w_{\Lambda,\hbar}(t\tau_{\hbar}),a\ra:=\left\la v_{\hbar}(t\tau_{\hbar}),\Oph^w\left(a\left(x,\xi,\frac{H_{\Lambda}(x,\xi)}{\eps_{\hbar}}
\right)\right)v_{\hbar}(t\tau_{\hbar})\right\ra.$$
The symbol involved belongs to the class of symbols $S^{0,0}_{\text{per}}$ amenable to pseudodifferential calculus on $\IT^2$. It will be useful to keep in mind 
Remark~\ref{r:semiclas} throughout this section.

\begin{rema}\label{r:splitting} We emphasize that, for $a$ in $\ml{C}^{\infty}_c(T^*\IT^2)$, one has
 $$\la w_{\hbar}(t\tau_{\hbar}),a\ra=\la w_{\Lambda,\hbar}(t\tau_{\hbar}),a\ra.$$
\end{rema}

Our first step is to explain how to extract converging subsequences following more or less standard procedures~\cite{Ge91, Ma09, AM10, Zw12}. 
For the sake of completeness, 
we briefly recall it. For that purpose, we denote by
$$\mathcal{B}:=\mathcal{C}^0_0(\IR^2\times\widehat{\IR},\mathcal{C}^3(\IT^2)),$$
the space of continuous function on $\IR^2\times\widehat{\IR}$ with values in $\mathcal{C}^3(\IT^2)$ and which tends to $0$ at infinity. We endow this space with its natural topology of Banach space. According to Theorem~\ref{t:cald-vail}, one knows that, for every $a$ in $\ml{C}^{\infty}_c(\IR\times T^*\IT^2\times\widehat{\IR})$, one has
\begin{equation}\label{e:measure}|\la w_{\Lambda,\hbar}(t\tau_{\hbar}),a(t)\ra|\leq C\|a(t)\|_{\mathcal{B}}.\end{equation}
Thus, the map $t\mapsto w_{\Lambda,\hbar}(t\tau_{\hbar})$ defines a bounded sequence in $L^1(\IR,\ml{B})'$, and, after extracting a subsequence,
 one finds that there exists $\mu_{\Lambda}$ in $L^1(\IR,\ml{B})'$ such that, for every $a$ in $\ml{C}^{\infty}_c(\IR\times T^*\IT^2\times\widehat{\IR})$, one has
$$\lim_{\hbar\rightarrow 0^+}\int_{\IR\times T^*\IT^2\times\widehat{\IR}}a(t,x,\xi,\eta) w_{\Lambda,\hbar}(t\tau_{\hbar},dx,d\xi,d\eta)dt
=\int_{\IR\times T^*\IT^2\times\widehat{\IR}}a(t,x,\xi,\eta) \mu_{\Lambda}(dt,dx,d\xi,d\eta).$$
Thanks to~\eqref{e:measure}, recall that, for every $\theta$ in $\ml{C}^{\infty}_c(\IR)$ and for every $a$ in 
$\ml{C}^{\infty}_c(T^*\IT^2\times\widehat{\IR})$, one has
$$\left|\int_{\IR\times T^*\IT^2\times\widehat{\IR}}\theta(t)a(x,\xi,\eta) \mu_{\Lambda}(dt,dx,d\xi,d\eta)\right|\leq C\|\theta\|_{L^1}\|a\|_{\ml{B}}.$$
Hence, $\mu_{\Lambda}$ is absolutely continuous with respect to the $t$ variable, i.e. for every $\theta$ in $L^1(\IR)$ and 
every $a$ in $\ml{C}^{\infty}_c(T^*\IT^2\times\widehat{\IR})$, one has
$$\lim_{\hbar\rightarrow 0^+}\int_{\IR}\theta(t)\la w_{\Lambda,\hbar}(t\tau_{\hbar}),a\ra dt =\int_{\IR}\theta(t)\la\mu_{\Lambda}(t),a\ra dt,$$
where, for a.e. $t$ in $\IR$, $\mu_{\Lambda}(t)\in\ml{B}'$. 

\subsection{Proof of Proposition~\ref{p:first-properties}}

We will first prove that the linear functionals $\mu_\Lambda$ are positive.
To see this, take $a\in\ml{C}^\infty_c(T^*\IT^2\times \widehat{\IR})$ such that $a\geq 0$. Using G\aa{}rding inequality (Th.~4.32 in~\cite{Zw12}), we deduce that
$$\la w_{\Lambda,\hbar}(t\tau_{\hbar}),a\ra\geq \ml{O}(\hbar\eps_{\hbar}^{-1})=o(1);$$
\begin{rema}
 Note that the proof of the G\aa{}rding inequality in~\cite{Zw12} is given in the case of $\IR^d$. The extension to compact manifolds usually requires to deal 
 with symbols that decay in $\xi$ as we differentiate with respect to $\xi$. Yet, in the case of the torus, we can verify that this property remains true 
 for an observable $a$ all of whose derivatives are bounded (i.e. not necessarily decaying in $\xi$) as in $\IR^d$. For that purpose, one can start from the 
 G\aa{}rding inequality on $\IR^d$ and apply the arguments of the proof of~\cite[Th.~5.5]{Zw12} which shows $L^2$-boundedness of pseudodifferential of order $0$ on $\IT^d$.
\end{rema}

After integrating against a test function $\theta$ in $L^1(\IR)$ and passing to the limit $\hbar\rightarrow 0$, 
one finds that, for a.e. $t$ in $\IR$,$$\la \mu_{\Lambda}(t),a\ra\geq 0.$$
Let $\chi\in\ml{C}^\infty_c(\IR)$ be a smooth cutoff function, with values in $[0,1]$ which is equal to $1$ in a neighborhood of $0$. For every $a\in\ml{C}^\infty_c(T^*\IT^2\times \widehat{\IR})$, we write
$$a_R(x,\xi,\eta):=a(x,\xi,\eta)\chi\left(\frac{\eta}{R}\right).$$
We define
$$\la \tilde{\mu}_\Lambda (t), a\ra := \lim_{R\to\infty} \la \mu_\Lambda (t), a_R\ra.$$
This limit clearly exists if $a\geq 0$, since $a_R$ is increasing and $\mu_\Lambda (t)$ is positive a.e.. The existence in the general case follows from the fact that, in general, one always can write $a=a_1-a_2$ for some non-negative $a_{1},a_2\in\ml{C}^\infty_c(T^*\IT^2\times \widehat{\IR})$. Note that, by definition, $\tilde{\mu}_\Lambda (t)$ is a positive functional, and $\tilde{\mu}_\Lambda (t)\leq \mu_\Lambda (t)$ for a.e. $t\in\IR$. This implies that the functional
$$\ml{C}^\infty_c(T^*\IT^2\times \IR)\ni a\longmapsto \la \tilde{\mu}_\Lambda (t), a\ra\in\IC$$
is a positive distribution, and therefore extends to a positive (finite) Radon measure on $T^*\IT^2\times \IR$.
Finally, note that, for every $a\in \ml{C}^\infty_c(T^*\IT^2\times \widehat{\IR})$ that vanishes at $\eta = \pm\infty$, one has
\begin{equation}\label{e:compactdef}
\la \tilde{\mu}_\Lambda (t), a\ra = \la \mu_\Lambda (t), a\ra , \quad\text{ for a.e. } t\in\IR,
\end{equation}
since in this case $a_R$ converges to $a$  in $\ml{B}$ as $R\to +\infty$. This in particular shows that the positive functional
$$\tilde{\mu}^\Lambda (t) := \mu_\Lambda (t)- \tilde{\mu}_\Lambda (t),$$
verifies that $\la \tilde{\mu}^\Lambda (t), a\ra$ only depends on the values of $a$ at $\eta = \pm \infty$. This means that if the restriction $a_1,a_2\in\ml{C}^\infty_c(T^*\IT^2\times \widehat{\IR})$ to $\eta = \pm\infty$ coincide then $\la \tilde{\mu}^\Lambda (t), a_1\ra=\la \tilde{\mu}^\Lambda (t), a_2\ra$ for a.e. $t\in\IR$.  This implies the existence, for a.e. $t\in\IR$, of distributions $\tilde{\mu}^\Lambda_{\pm}(t)\in\ml{D}'(T^*\IT^2)$ such that:
$$\tilde{\mu}^\Lambda=\tilde{\mu}^\Lambda_+\otimes\delta_{+\infty}+\tilde{\mu}^\Lambda_-\otimes\delta_{-\infty}.$$
Finally, $\tilde{\mu}^\Lambda_{\pm}(t)$ are necessarily positive since $\tilde{\mu}^\Lambda(t)$ is. Therefore they can be extended to positive Radon measures. This concludes the proof that $\mu_\Lambda$ is a positive, finite Radon measure on $T^*\IT^2\times \widehat{\IR}$ and one checks that $\tilde{\mu}_\Lambda (t)=\mu_\Lambda (t)\rceil_{T^*\IT^2\times \IR}$ and that $\tilde{\mu}^\Lambda(t)=\mu_\Lambda (t)\rceil_{T^*\IT^2\times \{\pm \infty\}}$. Thanks to the frequency assumption~\eqref{e:shosc}, one has, for a.e. $t$ in $\IR$,
\begin{equation}\label{e:zero-section-2-micro}
\mu_{\Lambda}(t)(\{\xi=0\})=0. 
\end{equation}

\begin{rema}\label{r:proj} Remark~\ref{r:splitting} implies that, for a.e. $t$ in $\IR$, the time-dependent semiclassical measure $\mu(t)$ can be obtained by 
\begin{equation}\label{e:splitting}
\mu(t)=\int_{\widehat{\IR}}\mu_{\Lambda}(t,\cdot,d\eta). 
\end{equation}
\end{rema} 
 
Concerning the support of $\tilde{\mu}_{\Lambda}(t)$, we let $a$ be an element in $\ml{C}^{\infty}_c(T^*\IT^2\times\IR)$ 
whose support does not intersect $\IT^2\times\Lambda^{\perp}\times\IR$. Using Remark \ref{r:semiclas}:
$$\Oph^w\left(a\left(x,\xi,\frac{H_{\Lambda}(x,\xi)}{\eps_{\hbar}}\right)\right)=\Op^w_{\hbar\eps_{\hbar}^{-1}}
\left(a\left(x,\eps_{\hbar}\xi,H_{\Lambda}(x,\xi)\right)\right).$$
Hence, this operator is equal to $0$ when $\hbar$ is small enough (thanks to our assumption on the support of $a$). This concludes the proof of the first part 
of Proposition~\ref{p:first-properties}.

Let us now discuss invariance by the geodesic flow. Again, we start with the ``compact'' part and we fix $a$ to be an element 
in $\ml{C}^{\infty}_c(T^*\IT^2\times\IR)$. Using composition rules for pseudodifferential operators, we write
\begin{align*}
 \frac{d}{dt}\la w_{\Lambda,\hbar}(t\tau_{\hbar}),a\ra  & =  \tau_{\hbar}\la w_{\Lambda,\hbar}(t\tau_{\hbar}),\xi.\partial_xa\ra\\
 &  +\frac{i\tau_{\hbar}\eps_{\hbar}^2}{\hbar}
\la v_{\hbar}(t\tau_{\hbar}), \left[V,\Op^w_{\hbar\eps_{\hbar}^{-1}}\left(a(x,\eps_{\hbar}\xi,H_{\Lambda}(x,\xi))\right)\right]v_{\hbar}(t\tau_{\hbar})\ra .
\end{align*}
Using Theorem~\ref{t:composition} (more specifically Remark~\ref{r:symmetry}) one more time, we have that
$$\left[V,\Op^w_{\hbar\eps_{\hbar}^{-1}}\left(a(x,\eps_{\hbar}\xi,H_{\Lambda}(x,\xi))\right)\right]
=-\frac{\hbar}{i\eps_{\hbar}}\Oph^w\left(\frac{\mathfrak{e}_{\Lambda}}{L_{\Lambda}}.
\partial_xV \partial_{\eta}a\left(x,\xi,\frac{H_{\Lambda}(x,\xi)}{\eps_{\hbar}}\right)\right)
+\ml{O}(\hbar(1+\hbar^2(\eps_{\hbar})^{-3})).$$
Combining these two identities to the fact $\hbar\eps_{\hbar}^{-1}=o(1)$, we find that
$$
 \frac{d}{dt}\la w_{\Lambda,\hbar}(t\tau_{\hbar}),a\ra  = \tau_{\hbar}\left(\left\la w_{\Lambda,\hbar}(t\tau_{\hbar}),\xi.\partial_xa
-\eps_{\hbar}\frac{\mathfrak{e}_{\Lambda}}{L_{\Lambda}}.\partial_xV \partial_{\eta}a \right\ra+o(\hbar)+\ml{O}(\eps_{\hbar}^2)\right).$$
Let now $\theta$ be an element in $\ml{C}^1_c(\IR)$. Integrating the previous equality against $\theta$ and integrating by parts, we find
$$\int_{\IR}\theta(t)\left\la w_{\Lambda,\hbar}(t\tau_{\hbar}),\xi.\partial_xa
-\eps_{\hbar}\frac{\mathfrak{e}_{\Lambda}}{L_{\Lambda}}.\partial_xV \partial_{\eta}a \right\ra dt=\ml{O}(\tau_{\hbar}^{-1})+o(\hbar)+\ml{O}(\eps_{\hbar}^2),$$
which implies the result for every $a$ in $\ml{C}^{\infty}_c(T^*\IT^2\times\IR)$ when we let $\hbar$ goes to $0$. Note that we used the first part of the 
Calder\'on-Vaillancourt Theorem~\ref{t:cald-vail} to bound the $\eps_{\hbar}$ term on the left hand side of this equality.

It now remains to treat the part at infinity. Let $a$ be an element in $\ml{C}^{\infty}_c(\IT^2\times\widehat{\IR})$. For every $R\geq 1$, we set
$$a^{R}(x,\xi,\eta):=a(x,\xi,\eta)\left(1-\chi\left(\frac{\eta}{R}\right)\right).$$The same argument as before allows to prove that, 
for every $\theta$ in $\ml{C}^1(\IR)$, one has
$$\int_{\IR}\theta(t)\left\la w_{\Lambda,\hbar}(t\tau_{\hbar}),(\xi.\partial_xa)^R
-\eps_{\hbar}\frac{\mathfrak{e}_{\Lambda}}{L_{\Lambda}}.\partial_xV \partial_{\eta}a^R \right\ra dt=o(1).$$
Thus, we can take the limit $\hbar\rightarrow 0$ and conclude the proof by letting $R$ goes to $+\infty$.

\subsection{Invariance and propagation of $2$-microlocal distributions} We now turn to the proofs of our main statements, namely 
Theorems~\ref{t:propag-inv} and~\ref{t:maintheo-infty}. Recall that a key ingredient of our proof in the Zoll case was an averaging argument of Weinstein~\cite{We77}. 
Here, it will be transposed by defining the differential operators
$$D_{\Lambda}:=\frac{1}{i}\frac{\mathfrak{e}_{\Lambda}}{L_{\Lambda}}.\nabla\ \ \text{and}\ \ D_{\Lambda}^{\perp}:=\frac{1}{i}\frac{\mathfrak{e}_{\Lambda}^{\perp}}{L_{\Lambda}}.\nabla$$
associated with the Hamiltonians $H_{\Lambda}$ and $H_{\Lambda}^{\perp}$. Clearly
\begin{equation}\label{e:laplacian}
 -\Delta=(D_{\Lambda}^{\perp})^2+D_{\Lambda}^2.
\end{equation}
Recall also that, for every smooth compactly supported function $b$ on $T^*\IT^2$, 
the Egorov theorem is exact for these operators and it tells us that
\begin{equation}\label{e:weinstein}\Oph^w(\ml{I}_{\Lambda}(b))=
\frac{1}{ L_{\Lambda}}\int_0^{ L_{\Lambda}} e^{isD_{\Lambda}^{\perp}}\Oph^w(b)e^{-isD_{\Lambda}^{\perp}}ds.\end{equation}
and that
\begin{equation}\label{e:commute}[D_{\Lambda}^{\perp},\Oph^w(\ml{I}_{\Lambda}(b))]=0,\end{equation}
which is at the heart of Weinstein's averaging method. 

\subsubsection{Proof of Theorem~\ref{t:propag-inv}} Let $a$ be an element in $\ml{C}^{\infty}_c(T^*\IT^2\times\IR)$. We start our proof 
by computing the derivative of the $2$-microlocal Wigner distribution. One has
$$\frac{d}{dt}\la w_{\Lambda, \hbar}(t\tau_{\hbar}),\ml{I}_{\Lambda}(a)\ra=\frac{i\tau_{\hbar}}{\hbar}\left\la v_{\hbar}(t\tau_{\hbar}),\left[\frac{\hbar^2}{2}(D_{\Lambda}^{\perp})^2+\frac{\hbar^2}{2}D_{\Lambda}^2+\eps_{\hbar}^2V,\Oph^w\left(a_{\Lambda,\hbar}\right)\right]v_{\hbar}(t\tau_{\hbar})\right\ra,$$
where
$$a_{\Lambda,\hbar}(x,\xi):=\ml{I}_{\Lambda}(a)\left(x,\xi,\frac{H_{\Lambda}(x,\xi)}{\eps_{\hbar}}\right).$$
Using~\eqref{e:commute}, we deduce that
$$\frac{d}{dt}\la w_{\Lambda, \hbar}(t\tau_{\hbar}),\ml{I}_{\Lambda}(a)\ra=\frac{i\tau_{\hbar}}{\hbar}\left\la v_{\hbar}(t\tau_{\hbar}),\left[\frac{\hbar^2}{2}D_{\Lambda}^2+\eps_{\hbar}^2V,\Oph^w\left(a_{\Lambda,\hbar}\right)\right]v_{\hbar}(t\tau_{\hbar})\right\ra,$$
Thanks to the commutation properties of the Weyl quantization from Remark~\ref{r:symmetry}, one has
$$\frac{d}{dt}\la w_{\Lambda, \hbar}(t\tau_{\hbar}),\ml{I}_{\Lambda}(a)\ra=\ml{O}(\eps_{\hbar}\tau_{\hbar}(\eps_{\hbar}+\hbar^2\eps_{\hbar}^{-2}))$$
\begin{equation}\label{e:deriv-2-micro}+\eps_{\hbar}\tau_{\hbar}\left\la v_{\hbar}(t\tau_{\hbar}),\Oph^w\left(\frac{H_{\Lambda}(x,\xi)}{\eps_{\hbar}}\frac{\mathfrak{e}_{\Lambda}.\partial_x\ml{I}_{\Lambda}(a)(x,\xi,H_{\Lambda}(x,\xi)/\eps_{\hbar})}{L_{\Lambda}}-\partial_{\eta}\ml{I}_{\Lambda}(a)\frac{\mathfrak{e}_{\Lambda}.\partial_xV}{L_{\Lambda}}\right)v_{\hbar}(t\tau_{\hbar})\right\ra.\end{equation}
Our assumption on the size of the perturbation ($\eps_{\hbar}\gg\hbar$) ensures that the remainder is in fact of order $o(\eps_{\hbar}\tau_{\hbar})$.

We now distinguish three regimes.

First, we suppose that $\eps_{\hbar}\tau_{\hbar}\rightarrow 0$ as $\hbar\rightarrow 0^+$. Thanks to the Calder\'on-Vaillancourt 
Theorem~\ref{t:cald-vail}, we can verify that the right hand-side of equality~\eqref{e:deriv-2-micro} 
is in fact $o(1)$ uniformly for $t$ in $\IR$. Letting $\hbar\rightarrow 0$, one finds that, for a.e. $t$ in $\IR$,
$$\mu_{\Lambda}(t)(\ml{I}_{\Lambda}(a))=\mu_{\Lambda}^0(\ml{I}_{\Lambda}(a)).$$
Combining Proposition~\ref{p:first-properties} with~\eqref{e:zero-section-2-micro}, one has then 
$\mu_{\Lambda}(t)(a)=\mu_{\Lambda}^0(\ml{I}_{\Lambda}(a))$ for a.e. $t$ in $\IR$, which proves point (1) of the Theorem.

Suppose now that $\tau_{\hbar}\eps_{\hbar}=1$. Letting $\hbar\rightarrow 0$, the limit measure satisfies the following transport equation, for all $\theta\in\ml{C}^1_c(\IR)$:
$$-\int_{\IR}\theta'(t)\mu_{\Lambda}(t)(\ml{I}_{\Lambda}(a))dt=\int_{\IR}\theta(t)
\mu_{\Lambda}(t)\left(\eta\frac{\mathfrak{e}_{\Lambda}.\partial_x\ml{I}_{\Lambda}(a)}{L_{\Lambda}}
-\partial_{\eta}\ml{I}_{\Lambda}(a)\frac{\mathfrak{e}_{\Lambda}.\partial_xV}{L_{\Lambda}}\right)dt.$$
Using again Proposition~\ref{p:first-properties} with~\eqref{e:zero-section-2-micro}, one deduces that
$$\partial_t\mu_{\Lambda}(t)(\ml{I}_{\Lambda}(a))=\mu_{\Lambda}(t)\left(\eta\frac{\mathfrak{e}_{\Lambda}.\partial_x\ml{I}_{\Lambda}(a)}{L_{\Lambda}}-\partial_{\eta}\ml{I}_{\Lambda}(a)\frac{\mathfrak{e}_{\Lambda}.\partial_x\ml{I}_{\Lambda}(V)}{L_{\Lambda}}\right).$$
This proves point (2) of the Theorem.

Finally, we suppose that $\tau_{\hbar}\eps_{\hbar}\rightarrow+\infty$. Let $\theta$ be an element 
in $\ml{C}^{1}_c(\IR)$. We integrate one more time equality~\eqref{e:deriv-2-micro} against $\theta$, and 
we make an integration by parts on the left-hand side of the equality. Then, we make use of the Calder\'on-Vaillancourt 
Theorem~\ref{t:cald-vail} to bound the left-hand-side. After letting $\hbar$ goes to $0$, one finds that, for every $\theta$ in $\ml{C}^1_c(\IR)$, 
$$\int_{\IR}\theta(t)\mu_{\Lambda}(t)\left(\eta\frac{\mathfrak{e}_{\Lambda}.\partial_x\ml{I}_{\Lambda}(a)}{L_{\Lambda}}-\partial_{\eta}\ml{I}_{\Lambda}(a)\frac{\mathfrak{e}_{\Lambda}.\partial_x\ml{I}_{\Lambda}(V)}{L_{\Lambda}}\right)dt=0,$$
where we used one more time Proposition~\ref{p:first-properties} with~\eqref{e:zero-section-2-micro} in order to replace 
$V$ by its $\Lambda$-average $\ml{I}_{\Lambda}(V)$. This implies point (3) of the Theorem.

\subsubsection{Proof of Theorem~\ref{t:maintheo-infty}}

Let now $a$ be an element in $\ml{C}^{\infty}_c(\IR^2\times\widehat{\IR})$ and let $k$ be an element in $\Lambda-\{0\}$. We fix 
$\chi_1(\eta)$ to be a smooth function on $\IR$ which is equal to $1$ for $\eta\geq 1$ and to $0$ for $\eta\leq 1/2$. For every $R\geq 1$, we set
$$a^{R,k}_{\pm}(x,\xi,\eta):=e^{-2i\pi k.x}a(\xi,\eta)\chi_1\left(\pm\frac{\eta}{R}\right).$$
\begin{rema}\label{r:IPP} Let $\theta$ be an element in $\ml{C}^1_c(\IR)$. One has
 $$\int_{\IR}\theta(t)\frac{d}{dt}\left\la w_{\Lambda,\hbar}(t\tau_{\hbar}),\frac{1}{\eta}a_{\pm}^{R,k}\right\ra dt=
 -\int_{\IR}\theta'(t)\left\la w_{\Lambda,\hbar}(t\tau_{\hbar}),\frac{1}{\eta}a_{\pm}^{R,k}\right\ra dt.$$
Thanks to the Calder\'on-Vaillancourt Theorem~\ref{t:cald-vail}, one knows that
$$\left\|\Oph^w\left(\chi\left(\frac{H_{\Lambda}(x,\xi)}{R\eps_{\hbar}}\right)
a\left(\xi,\frac{H_{\Lambda}(x,\xi)}{\eps_{\hbar}}\right)e^{-2 i\pi k.x}\frac{\eps_{\hbar}}{H_{\Lambda}(x,\xi)}\right)
\right\|_{L^2\rightarrow L^2}=\ml{O}(R^{-1}).$$
Thus, one has
$$\int_{\IR}\theta(t)\frac{d}{dt}\left\la w_{\Lambda,\hbar}(t\tau_{\hbar}),\frac{1}{\eta}a_{\pm}^{R,k}\right\ra dt=\ml{O}(R^{-1}).$$
\end{rema}

In order to prove the proposition, we will now compute explicitely the derivative of 
$\left\la w_{\Lambda,\hbar}(t\tau_{\hbar}),\frac{1}{\eta}a_{\pm}^{R,k}\right\ra.$ For that purpose, we need to compute the following bracket:
$$\left[-\frac{\hbar^2\Delta}{2}+\eps_{\hbar}^2V,\Oph^w\left(a_{\pm}^{R,k}\left(x,\xi,\frac{H_{\Lambda}(x,\xi)}{\eps_{\hbar}}\right)
\frac{\eps_{\hbar}}{H_{\Lambda}(x,\xi)}\right)\right].$$
Using again~\eqref{e:commute}, this commutator is in fact equal to
$$\left[\frac{\hbar^2D_{\Lambda}^2}{2}+\eps_{\hbar}^2V,\Oph^w\left(a_{\pm}^{R,k}\left(x,\xi,\frac{H_{\Lambda}(x,\xi)}{\eps_{\hbar}}
\right)\frac{\eps_{\hbar}}{H_{\Lambda}(x,\xi)}\right)\right].$$
We split this commutator in two parts. Thanks to remark~\ref{r:symmetry}, one has
$$\left[\frac{\hbar^2D_{\Lambda}^2}{2},\Oph^w\left(a_{\pm}^{R,k}\left(x,\xi,\frac{H_{\Lambda}(x,\xi)}{\eps_{\hbar}}
\right)\frac{\eps_{\hbar}}{H_{\Lambda}(x,\xi)}\right)\right]
=-2\pi\hbar\eps_{\hbar}\Oph^w\left(\frac{\mathfrak{e}_{\Lambda}}{L_{\Lambda}}.k a_{\pm}^{R,k}\left(x,\xi,\frac{H_{\Lambda}(x,\xi)}{\eps_{\hbar}}\right)\right).$$
For the other part of the commutator, we use one more time the commutation rule for pseudodifferential operators and the Calder\'on Vaillancourt Theorem~\ref{t:cald-vail}. We find that
$$\left[V,\Oph^w\left(a_{\pm}^{R,k}\left(x,\xi,\frac{H_{\Lambda}(x,\xi)}{\eps_{\hbar}}
\right)\frac{\eps_{\hbar}}{H_{\Lambda}(x,\xi)}\right)\right]=
\ml{O}_{L^2\rightarrow L^2}\left(\hbar\eps_{\hbar}^{-1}R^{-1}+\hbar+\hbar^2\eps_{\hbar}^{-2}\right).$$
As $\hbar\eps_{\hbar}^{-1}\rightarrow 0$, we finally get that
$$\frac{d}{dt}\left\la w_{\Lambda,\hbar}(t\tau_{\hbar}),\frac{1}{\eta}a_{\pm}^{R,k}\right\ra=
-\frac{2\pi\tau_{\hbar}\eps_{\hbar}\mathfrak{e}_{\Lambda}.k}{L_{\Lambda}}\left\la w_{\Lambda,\hbar}(t\tau_{\hbar}),a_{\pm}^{R,k}\right\ra+
\ml{O}(\tau_{\hbar}\eps_{\hbar} R^{-1})+o(\tau_{\hbar}\eps_{\hbar}).$$
Let now $\theta$ be an element in $\ml{C}^1_c(\IR)$. We integrate these expressions against $\theta$. Using Remark~\eqref{r:IPP} 
and making the assumption that $\tau_{\hbar}\eps_{\hbar}\geq 1$, we obtain
$$\forall k\in\Lambda-\{0\},\ \int_{\IR}\theta(t)\left\la w_{\Lambda,\hbar}(t\tau_{\hbar}),a_{\pm}^{R,k}\right\ra dt=o(1)+\ml{O}(R^{-1}).$$
We now let $\hbar$ goes to $0$, and we get that, for every $R>0$,
$$\forall k\in\Lambda-\{0\},\ \int_{\IR}\theta(t)\left\la \mu_{\Lambda}(t),a_{\pm}^{R,k}\right\ra dt=\ml{O}(R^{-1}).$$
To get the conclusion, we let $R$ goes to $+\infty$.

\begin{rema}\label{r:constant} From this Theorem, we deduce that, for every 
$a(x,\xi,\eta)$ in $\mathcal{C}^{\infty}_c(T^*\IT^2\times\widehat{\IR})$ and for a.e. $t$ in $\IR$,
$$\tilde{\mu}^{\Lambda}(t)(\ml{I}_{\Lambda}(a))=\int_{T^*\IT^2\times\{\pm\infty\}}\widehat{a}_0(\xi,\eta)\mu_{\Lambda}(t, d\xi,d\eta).$$ 
\end{rema}

\appendix

\section{Regularity of bi-invariant measures}
\label{a:biinvariant}

In this appendix, we fix $\Lambda$ a primitive sublattice of $\IZ^2$ of rank $1$, and we aim at 
analyzing the regularity of the set of finite measures on $T^*\IT^2$ which are invariant 
by the Hamiltonian flows\footnote{By making a slight abuse of notation, we shall identify $\varphi_{p_\Lambda}^t$, a flow \textit{a priori} defined on $\IT^2\times\Lambda^\bot\times\IR$, to a flow on $T^*\IT^2$ via the diffeomorphism \eqref{e:diffeo}. Recall that $\varphi_{H_{\Lambda}^{\perp}}^t$ and $\varphi_{p_{\Lambda}}^t$ commute.} $\varphi_{H_{\Lambda}^{\perp}}^t$ and $\varphi_{p_{\Lambda}}^t$. 
We will now recall the results from section~4 of~\cite{MaRi15} and explain how they can be adapted to the present framework. We refer the reader 
to this reference for the detailed proofs. We introduce the critical set in the direction of $\Lambda$:
$$\text{Crit}_{\Lambda}(V):=\{(x,\xi)\in T^*\IT^2: H_{\Lambda}(x,\xi)=0\ \text{and}\ \partial_x\ml{I}_{\Lambda} (V)=0\}.$$
This is a closed subset of $T^*\IT^2$ which is invariant by the Hamiltonian flows $\varphi_{H_{\Lambda}^{\perp}}^t$ and $\varphi_{p_{\Lambda}}^t$, and we introduce its complement
$$\ml{R}(\Lambda):=T^*\IT^2-\text{Crit}_{\Lambda}(V).$$
The map
$$\phi:\IR^2 \times  \ml{R}(\Lambda) \ni (s,t,x,\xi) \longmapsto \varphi_{H_{\Lambda}^{\perp}}^s \circ \varphi_{p_{\Lambda}}^t(x,\xi) \in \ml{R}(\Lambda),$$
is a group action of $\IR^2$ on $\ml{R}(\Lambda)$. Moreover, for any $(x_0,\xi_0)\in \ml{R}(\Lambda)$, the map
$$\phi_{x_0,\xi_0}:\IR^2\ni (s,t) \longmapsto \varphi_{H_{\Lambda}^{\perp}}^s \circ \varphi_{p_{\Lambda}}^t(x_0,\xi_0) \in \ml{R}(\Lambda), $$
is an immersion. Therefore, the stabilizer group $G_{x_0,\xi_0}$ of $(x_0,\xi_0)$ under $\phi$ is discrete. 
This proves that the orbits of the action $\phi$ are either diffeomorphic to the torus $\IT^2$, to the cylinder $\IT\times \IR$ or to $\IR^2$. On the 
other hand, the moment map,
$$\Phi : \ml{R}(\Lambda) \ni (x,\xi) \longmapsto (H_{\Lambda}^{\perp}(x,\xi), p_{\Lambda}(x,\xi)) \in \IR^2,$$
is a submersion, and, for every $(H,J)\in\Phi(\ml{R}(\Lambda))$ the level set 
$$\ml{L}_{(H,J)}:=\Phi^{-1}(H,J),$$
is a smooth submanifold of $\ml{R}(\Lambda)$ of dimension two. To summarize, the couple $(H_{\Lambda}^{\perp}, p_{\Lambda})$ 
forms a completely integrable system on $\ml{R}(\Lambda)$, and the map $\phi_{x_0,\xi_0}$ induces a diffeomorphism:
$$\forall (x_0,\xi_0)\in\ml{R}(\Lambda),\quad
\phi_{x_0,\xi_0}:\IR^2/G_{x_0,\xi_0}\longrightarrow \ml{L}_{(H_0,J_0)}^{x_0,\xi_0},\quad \text{ for } (H_0,J_0):=\Phi(x_0,\xi_0).$$ 
Here, $\ml{L}_{(H_0,J_0)}^{x_0,\xi_0}$ denotes the connected component of $\ml{L}_{(H_0,J_0)}$ that contains $(x_0,\xi_0)$. Therefore, if $\ml{L}_{(H_0,J_0)}^{x_0,\xi_0}$ is compact then it is an embedded Lagrangian torus in $T^*\IT^2$. In that case, we shall write $\IT_{x_0,\xi_0}^2:=\IR^2 / G_{x_0,\xi_0}$. In the following, we denote by $\ml{R}_c(\Lambda)$ the set formed by those $(x,\xi)\in\ml{R}(\Lambda)$ such that $\ml{L}_{\Phi(x,\xi)}^{x,\xi}$ is compact. Mimicking the proof of proposition~4.2 in~\cite{MaRi15}, one can show that the following holds:
\begin{prop}\label{p:biinv}
Let $\mu$ be a probability measure on $\ml{R}(\Lambda)$ that is invariant by $\varphi_{H_{\Lambda}^{\perp}}^t$ and $\varphi_{p_{\Lambda}}^t$. Set $\overline{\mu}:=\Phi_*\mu$. Then, for every $a\in\ml{C}_c(\ml{R}(\Lambda))$, one has
$$ \int_{\ml{R}(\Lambda)}a(x,\xi)\mu(dx,d\xi)=\int_{\Phi(\ml{R}(\Lambda))}\int_{\ml{L}_{(H,J)}}a(x,\xi)\lambda_{H,J}(dx,d\xi)\overline{\mu}(dH,dJ),$$
where, for $(H,J)\in\Phi(\ml{R}(\Lambda))$, the measure $\lambda_{H,J}$ is a convex combination of the (normalized) Haar measures on the tori $\ml{L}_{(H,J)}^{x_0,\xi_0}$ for $(x_0,\xi_0)\in\ml{L}_{(H,J)}\cap\ml{R}_c(\Lambda)$. In particular, for every $(x,\xi)$ in $\ml{R}(\Lambda)$, one has
$$\mu\left(\left\{\varphi_{H_{\Lambda}^{\perp}}^s(x,\xi): 0\leq s\leq L_{\Lambda}\right\}\right)=0.$$
\end{prop}
An explicit formula for the restriction of the measure $\lambda_{H,J}$ to a connected component $\ml{L}_{(H,J)}^{x,\xi}$ with $(x,\xi)\in\ml{R}_c(\Lambda)\cap\ml{L}_{(H,J)}$ is the following:
\begin{equation}\label{e:haar}
\int_{\ml{L}_{(H,J)}^{x_0,\xi_0}}a(x,\xi)\lambda_{H,J}(dx,d\xi)=c\int_{\IT_{x_0,\xi_0}^2}a(\phi_\rho(s,t))dsdt,
\end{equation}
for some constant $c\in [0,1]$.

We will now discuss the regularity of the projections of bi-invariant measures following the proof from paragraph~4.2 in~\cite{MaRi15}. 
We denote by $\Pi:T^*\IT^2\rightarrow \IT^2$ the canonical projection. The main result from section~4 in~\cite{MaRi15} was the following

\begin{theo}\label{t:abscon}
Let $\mu$ be a probability measure on $\ml{R}(\Lambda)$ that is invariant by $\varphi_{H_{\Lambda}^{\perp}}^t$ and $\varphi_{p_{\Lambda}}^t$. Then, 
$\nu:=\Pi_*\mu$ is a probability measure on $\IT^2$ that is absolutely continuous with respect to the Lebesgue measure. 
\end{theo}

Denote by $\ml{N}(\Lambda)$ the convex closure of the set of measures $\delta_{\Pi\circ\Gamma}$ where $\Gamma\subset T^*\IT^2$ 
ranges over the orbits of $\varphi_{H_{\Lambda}^{\perp}}$ that are contained in $\text{Crit}_{\Lambda}(V)$. A direct consequence of 
the previous Theorem is the following: 
\begin{coro}\label{c:struct}
The projection $\nu:=\Pi_*\mu$ of a probability measure $\mu$ on $T^*\IT^2$ that is invariant by $\varphi_{H_{\Lambda}^{\perp}}^t$ and $\varphi_{p_{\Lambda}}^t$ can be decomposed as:
$$\nu = f \operatorname{vol} + \alpha \nu_{\emph{sing}}$$
where $f\in L^1(\IT^2)$, $\alpha\in[0,1]$ and $\nu_{\emph{sing}}\in\ml{N}(\Lambda)$.
\end{coro}
Note that, for a ``generic'' choice of $V$, the set of points $x$ satisfying $\partial_x\ml{I}_{\Lambda}(V)=0$ consists of finitely many closed geodesics of $\IT^2$. In particular, $\nu_{\text{sing}}$ is a finite combination of measures carried by closed geodesics.
\begin{proof} As it is simple to explain in the current framework, we briefly explain how the proof of Theorem~4.6 in~\cite{MaRi15} can be adapted to prove Theorem~\ref{t:abscon} -- see also Lemma~2.1 in~\cite{BiPo89}. Recall that it is sufficient to fix some $(x_0,\xi_0)$ in $\ml{R}_{c}(\Lambda)$ and to prove that the set of points where
$$\phi_{x_0,\xi_0}:(s,t)\in\IT^2_{x_0,\xi}\mapsto \Pi\circ \varphi_{H_{\Lambda}^{\perp}}^s \circ \varphi_{p_{\Lambda}}^t(x_0,\xi_0)\in \IT^2$$
is not a local diffeomorphism is made of finitely many disjoint $\ml{C}^1$ closed curves. Such curves are called caustics. 
This can be proved as follows. One can verify that the points where we do not have a local diffeomorphism are defined by the points $(s,t)$ satisfying 
$$H_{\Lambda}\left(\phi_{x_0,\xi_0}(s,t)\right)=0.$$ 
Note that, for every $s$ in $\IR$,
$$H_{\Lambda}\left(\varphi_{p_{\Lambda}}^t(x_0,\xi_0)\right)=H_{\Lambda}\left(\phi_{x_0,\xi_0}(s,t)\right).$$
As $(x_0,\xi_0)$ belongs to the $\varphi_{p_{\Lambda}}^t$-invariant set $\ml{R}(\Lambda)$, we know that
$$\partial_x\ml{I}_{\Lambda}(V)\left(\varphi_{p_{\Lambda}}^t(x_0,\xi_0)\right)\neq 0.$$
Thus, from the Hamilton-Jacobi equations, we deduce that there exists a small open neighborhood $(t-\eta,t+\eta)$ of $t$ such that, for every $t'\in(t-\eta,t+\eta)-\{t\}$, 
$$H_{\Lambda}\circ\varphi_{p_{\Lambda}}^{t'}(x_0,\xi_0)\neq 0.$$
In particular, there are ony finitely many values of $t$ such that $H_{\Lambda}\circ\varphi_{p_{\Lambda}}^t(x_0,\xi_0)\neq 0$ and thus, 
there are only finitely many closed curves on $\IT_{x_0,\xi_0}^2$ where the map $\phi_{x_0,\xi_0}$ is not a local diffeomorphism.
\end{proof}

\section{Background on semiclassical analysis}
\label{a:sc-an}

In this appendix, we give a brief reminder on semiclassical analysis and we refer to~\cite{Zw12} (mainly Chapters $1$ to $5$) for a more detailed exposition. Given $\hbar>0$ 
and $a$ in $\ml{S}(\IR^{2d})$ (the Schwartz class), one can define the Weyl quantization of $a$ as follows:
$$\forall u\in\ml{S}(\IR^d),\ \Oph^w(a)u(x):=\frac{1}{(2\pi\hbar)^d}\iint_{\IR^{2d}}e^{\frac{i}{\hbar}\la x-y,\xi\ra}a\left(\frac{x+y}{2},\xi\right)u(y)dyd\xi.$$

This definition can be extended to any observable $a$ with uniformly bounded derivatives, i.e. such that for every $\alpha\in\IN^{2d}$, there exists $C_{\alpha}>0$ such that 
$\sup_{x,\xi}|\partial^{\alpha}a(x,\xi)|\leq C_{\alpha}$. More generally, we will use the convention, for every $m\in\IR$ and every $k\in\IZ$,
$$S^{m,k}:=\left\{(a_{\hbar}(x,\xi))_{0<\hbar\leq 1}:\ \forall (\alpha,\beta)\in\IN^d\times\IN^d,\ \sup_{(x,\xi)\in\IR^{2d}; 0<\hbar\leq 1}|\hbar^k\la\xi\ra^{-m}\partial_x^{\alpha}\partial_{\xi}^{\beta}a_{\hbar}(x,\xi)|<+\infty\right\},$$
where $\la\xi\ra:=(1+\|\xi\|^2)^{1/2}$.  
For such symbols, $\Oph^w(a)$ defines a continuous operator $\ml{S}(\IR^d)\rightarrow\ml{S}(\IR^d)$ which acts by duality on $\ml{S}'(\IR^d)$. 

\begin{rema} We also note that we have the following relation that we use at different stages of our proof:
\begin{equation}\label{e:change-variable}
\forall\delta>0,\ \forall a\in S^{m,k}, \Oph^w(a(x,\xi))=\Op^w_{\hbar\delta^{-1}}(a(x,\delta\xi)).
\end{equation}
\end{rema}

Among the above symbols, we distinguish the family of $\IZ^d$-periodic symbols that we denote by $S^{m,k}_{per}$. Note that any $a$ in $\ml{C}^{\infty}(T^*\IT^d)$ 
(with bounded derivatives) defines an element in $S^{0,0}_{per}$. Similarly to the proof of Th.~$4.19$ in~\cite{Zw12}, 
one can verify that, for any $a\in S^{m,k}_{per}$,
$$\Oph^w(a)(e_k)=\sum_{q\in\IZ^d}e_q\hat{a}_{q-k}(\pi\hbar(q+k)),$$
where $e_k(x):=e^{2i\pi k.x}$, and $\hat{a}_p(\xi):=\int_{\IT^d}a(x,\xi)e^{-2i\pi p.x}dx.$ In particular, for any $a\in S^{m,k}_{per}$, the operator $\Oph^w(a)$ maps trigonometric polynomials into a smooth $\IZ^d$-periodic function, and more generally any smooth $\IZ^d$-periodic function into a smooth $\IZ^d$-periodic function. Thus, for every $a$ in $S^{m,k}_{per}$, $\Oph^w(a)$ acts by duality on the space of distributions $\ml{D}'(\IT^d)$. An important feature of this quantization procedure is that it defines a bounded operator on $L^2(\IT^d)$:
\begin{theo}\label{t:cald-vail}[Calder\'on-Vaillancourt] There exists a constant $C_d>0$ and an integer $D>0$ such that, for every $a$ in $S^{0,0}_{per}$, one has, for every $0<\hbar\leq 1$,
$$\left\|\Oph^w(a)\right\|_{L^2(\IT^d)\rightarrow L^2(\IT^d)}\leq C_d\sum_{|\alpha|\leq d+1}\|\partial^{\alpha}_xa\|_{\infty},$$
and for every $a$ in $\ml{C}^{\infty}_c(T^*\IT^d)$,
$$\left\|\Oph^w(a)\right\|_{L^2(\IT^d)\rightarrow L^2(\IT^d)}\leq C_d\sum_{|\alpha|\leq D}\hbar^{\frac{|\alpha|}{2}}\|\partial^{\alpha}a\|_{\infty}.$$
\end{theo}
The second part of the Theorem follows from the fact that, when $a$ belongs to $\ml{C}^{\infty}_c(T^*\IT^d)$, $\Oph^w(a)$ defines a ``standard'' pseudodifferential operator on the manifold $\IT^d$. In particular, we can apply the usual Calder\'on-Vaillancourt Theorem (see e.g. Ch.~5 in~\cite{Zw12}) from which the second part of the Theorem follows. The advantage of the first part is that it allows to extend the Weyl quantization to more general symbols which may not vanish at infinity. Yet, this part is really specific to the case of the torus and we shall give a proof of it. 

\begin{proof}
The proof of the first part of the Theorem is an adaptation for the $\hbar$-Weyl quantization of the proof of Th.~4.8.1 in~\cite{Ru10} which was given for the $\hbar=1$-standard quantization. As we already observed it, we can write, for every trigonometric polynomial $u$ in 
$L^2(\IT^d)$, its Fourier decomposition $u=\sum_{k\in\IZ^d}\hat{u}_ke_k$, and one has then
$$\Op_{\hbar}^w\left(a\right)u=\sum_{k,q\in\IZ^d}\hat{u}_k\hat{a}_{q-k}(\pi\hbar(q+k))e_{q},$$
where $a(x,\xi)=\sum_{l\in\IZ^d}\hat{a}_l(\xi)e_l(x).$ Applying Plancherel equality, we get
$$\left\|\Op_{\hbar}^w\left(a\right)u\right\|_{L^2(\IT^d)}^2=\sum_{q\in\IZ^d}\left|\sum_{k\in\IZ^d}\hat{u}_k\hat{a}_{q-k}(\pi\hbar(q+k))\right|^2.$$
Thanks to Cauchy-Schwarz inequality, one has 
$$\left\|\Op_{\hbar}^w\left(a\right)u\right\|_{L^2(\IT^d)}^2\leq\sum_{q\in\IZ^d}\left(\sum_{k\in\IZ^d}|\hat{u}_k|^2|\hat{a}_{q-k}(\pi\hbar(q+k))|\right)
\left(\sum_{k'\in\IZ^d}|\hat{a}_{q-k'}(\pi\hbar(q+k'))|\right).$$
This implies that
$$\left\|\Op_{\hbar}^w\left(a\right)\right\|_{L^2(\IT^d)\rightarrow L^2(\IT^d)}^2\leq\sup_{q\in\IZ^d}\left(\sum_{k'\in\IZ^d}|\hat{a}_{q-k'}(\pi\hbar(q+k'))|\right)\times
 \sup_{k\in\IZ^d}\left(\sum_{q\in\IZ^d}|\hat{a}_{q-k}(\pi\hbar(q+k))|\right),$$
which concludes the proof of the lemma. 
\end{proof}

Another important feature of the Weyl quantization procedure is the composition formula:
\begin{theo}\label{t:composition}[Composition formula] Let $a\in S^{m_1,k_1}$ and $b\in S^{m_2,k_2}$. Then, one has, for any $0<\hbar\leq 1$
$$\Oph^w(a)\circ\Oph^w(b)=\Oph^w(a\sharp_{\hbar} b),$$
 in the sense of operators from $\ml{S}(\IR^d)\rightarrow\ml{S}(\IR^d)$, where $a\sharp_{\hbar} b$ has uniformly bounded derivatives, and, for every $N\geq 0$
$$a\sharp_{\hbar} b\sim \sum_{k=0}^N\frac{1}{k!}\left(\frac{i\hbar}{2}D\right)^k(a,b)+\ml{O}(\hbar^{N+1}),$$
where $D(a,b)(x,\xi)=(\partial_x\partial_{\nu}-\partial_{y}\partial_{\xi})(a(x,\xi)b(y,\nu))\rceil_{y=x,\nu=\xi}$.
\end{theo}
We refer to chapter~$4$ of~\cite{Zw12} for a detailed proof of this result. We observe that for $N=0$, the coefficient is given by the symbol $ab$, and for $N=1$,
it is given by $\frac{\hbar}{2i}\{a,b\}$, where $\{.,.\}$ is the Poisson bracket. As before, we can restrict this result to the case of periodic symbols, and we can check that the composition formula remains valid for operators acting on $\ml{C}^{\infty}(\IT^d)$.
\begin{rema}\label{r:symmetry} We note that the formula for the composed symbols is quite symmetric, and we have in fact the following useful property, for every $N\geq 0$,
 $$a\sharp_{\hbar} b-b\sharp_{\hbar} a\sim \sum_{k=0}^N\frac{2}{(2k+1)!}\left(\frac{i\hbar}{2}D\right)^{2k+1}(a,b)+\ml{O}(\hbar^{2N+3}),$$
Finally, note that, if $b(\xi)$ is a polynomial in $\xi$ of order $\leq 2$, one has, the exact formula:
$$a\sharp_{\hbar} b-b\sharp_{\hbar} a=\frac{\hbar}{2i}\{a,b\}.$$
\end{rema}

\end{document}